\documentclass[12pt]{article}
\usepackage{amsmath}

\usepackage[margin=1in, centering]{geometry}

\usepackage{placeins}
\usepackage{verbatim}

\usepackage{url}

\usepackage[dvips]{graphicx}
\usepackage{color}
\usepackage{amsfonts,amssymb,amsmath,}     
\usepackage[T1]{fontenc}
\usepackage{ae}

\newtheorem{theorem}{Theorem}

\newtheorem{coro}{Corollary}[theorem]

\newtheorem{lem}{Lemma}

\numberwithin{equation}{section} 
\numberwithin{table}{section} 
\numberwithin{figure}{section}

\renewcommand{\(}{\left(}
\renewcommand{\)}{\right)}

\allowdisplaybreaks[4]

\begin{document}
\title{Simultaneous Distribution of the Fractional Parts of Riemann Zeta Zeros}
\author{Kevin Ford, Xianchang Meng, and Alexandru Zaharescu}
\date{}
\maketitle
\begin{abstract}
In this paper, we investigate the simultaneous distribution of the fractional parts of $\{\alpha_1 \gamma, \alpha_2\gamma, \cdots, \alpha_n\gamma\}$, where $n\geq 2$,  $\alpha_1$, $\alpha_2$, $\cdots$, and $\alpha_n$ are fixed distinct positive real numbers and $\gamma$ runs over the imaginary parts of the non-trivial zeros of the Riemann zeta-function.

\end{abstract}

\let\thefootnote\relax
\let\thefootnote\relax

\section{Introduction and statement of results}
Let $\alpha$ be a fixed positive real number, and $\gamma$ run over the imaginary parts of the zeros of the Riemann zeta-function. We are interested in the distribution of the fractional parts $\{\alpha\gamma\}$. 
 Rademacher \cite{radema} was the first to consider this problem and he conjectured that, for a certain specific type of $\alpha$, 
there should be a "predominance of terms which fulfill $|\{\alpha\gamma\}-1/2|<1/4$".  
Since the fractional parts are uniformly distributed modulo $1$, as proved by Hlawka \cite{hla} in 1975, any
discrepancy must be very subtle.  
The first and third authors uncovered in \cite{fz} this delicate inequity in the fractional parts and 
 not only proved Rademacher correct, but also gave a much more precise measure for this phenomenon.

Let $\mathbb{T}=\mathbb{R}/\mathbb{Z}$ be the torus. Since the fractional parts $\{\alpha\gamma\}$ are uniformly distributed (mod 1) for any fixed $\alpha$, for all continuous functions $h: \mathbb{T}\rightarrow\mathbb{C}$, we have
$$\sum_{0<\gamma\leq T}h(\alpha\gamma)=N(T)\int_{\mathbb{T}}h(u)du+o(N(T)),$$
as $T\rightarrow\infty$, where $N(T)$ denotes the number of zeros $0<\gamma\leq T$.  We know that (\cite{Titch}, Theorem 9.4)
$$N(T)=\frac{T}{2\pi}\log\frac{T}{2\pi e}-\frac{T}{2\pi}+O(\log T),\qquad(T\geq 1).$$

The first and third authors \cite{fz} showed that, for a large class of functions $h: \mathbb{T}\rightarrow\mathbb{C}$, as $T\rightarrow\infty$,
$$\frac{1}{T}\sum_{0<\gamma\leq T}h(\alpha\gamma)-\frac{N(T)}{T}\int_{\mathbb{T}} h(u)du=\int_{\mathbb{T}} h(u)g_{\alpha}(u)du+o(1),$$
where $g_{\alpha}(u)$ is a function depending on the form of $\alpha$. From this result, we can see that the right hand side is close to a constant, and thus the discrepancy of the set $\{h(\alpha\gamma): 0<\gamma\leq T\}$ is of order $O(\frac{1}{\log T})$.
In \cite{fsz}, the first author, Soundararajan, and the third author established connections between the discrepancy of this set, Montgomery's pair correlation function and the distribution of primes in short intervals.

In the present paper, our goal is to generalize the results from \cite{fz} to the case of simultaneous distribution of fractional parts $(\{\alpha_1\gamma\}, \{\alpha_2\gamma\}, \cdots, \{\alpha_n \gamma\})$, where $\alpha_1,\ldots,\alpha_n$ are distinct, positive real numbers. 
As we will see below, a new phenomenon involving Diophantine approximation appears in the higher dimensional case $n\geq 2$. 

For  $n\geq 2$, consider a function $h(\mathbf{x})$ defined on $\mathbb{T}^n$, which has Fourier expansion
$$h(\mathbf{x})=\sum_{\mathbf{m}\in \mathbb{Z}^n} c_{\mathbf{m}} e^{2\pi i \mathbf{m}\cdot \mathbf{x}},$$
and assume that
\begin{equation}\label{n-dim-condt}
|c_{\mathbf{m}}|\ll \left(\frac{1}{\|\mathbf{m}\|+1}\right)^B,
\end{equation}
for some constant $B>n+2$, and $\|\cdot\|$ is the sup-norm.

Let $\boldsymbol{\alpha}\in \mathbb{R}^n$.   Suppose that there is
a full-rank $r\times n$ ($r\leq n$) matrix $M=(b_{ij})$ with integer entries,
 $\gcd(b_{i, 1}, b_{i, 2}, \cdots, b_{i, n})=1$ for all $1\leq i\leq r$, and
 such that
\begin{equation}\label{n-dim-condition}
M\boldsymbol{\alpha}^{\top}=P,
\end{equation}
where  $P=(\frac{a_1}{q_1}\frac{ \log p_1}{2\pi }, \cdots, \frac{a_r }{q_r}\frac{\log p_r}{2\pi })^{\top}$ for some integers $a_i$, $q_i$ and distinct primes $p_i$ ($i=1, \cdots, r$).  Among all possible such matrices, we take
one with maximal $r$,  in which case the above conditions ensure that $M$ 
is uniquely determined.  We then define
\begin{eqnarray}\label{n-dim-g1}
g_{\boldsymbol{\alpha}}(\mathbf{x})&=&-\frac{1}{\pi}\sum_{j=1}^r (\log p_j) \Re\left\{ \sum_{k\geq 1} p_j^{-\frac{a_j k}{2}} e^{-2kq_j \pi i(\mathbf{b_j}\cdot \mathbf{x})}  \right\}\nonumber\\
&=&-\frac{1}{\pi}\sum_{j=1}^r \frac{(\log p_j)(p_j^{\frac{a_j}{2}}\cos(2\pi q_j(\mathbf{b_j}\cdot \mathbf{x})-1)}{p_j^{a_j}-2p_j^{\frac{a_j}{2}}\cos(2\pi q_j(\mathbf{b_j}\cdot\mathbf{x})+1},
\end{eqnarray}
where $\mathbf{b_j}=(b_{j, 1}, b_{j, 2}, \cdots, b_{j, n})$ for $1\leq j\leq r$. 
If no such matrix exists, we simply set
\begin{equation}\label{n-dim-g2}
g_{\boldsymbol{\alpha}}(\mathbf{x})=0 \qquad (\mathbf{x} \in \mathbb{T}^n).
\end{equation}
Then we have the following.
\begin{theorem}\label{n-dim-aa}
	Let $n\geq 2$. Assume $h(\mathbf{x})$ satisfies condition (\ref{n-dim-condt}). Then,
	if $\boldsymbol{\alpha}\in \mathbb{R}^n$ satisfies 
	\begin{equation}\label{n-dim-alpha}
	|\mathbf{m}\cdot\boldsymbol{\alpha}|>\frac{C}{e^{\|\mathbf{m}\|}}
	\end{equation}
	for all $\mathbf{m}\in \mathbb{Z}^n$ and some fixed constant $C$,
	we have
	$$\lim_{T\rightarrow\infty}\frac{1}{T}\left(\sum_{0<\gamma\leq T}h(\gamma\boldsymbol{\alpha})-N(T)\int_{\mathbb{T}^n} h\right)=\int_{\mathbb{T}^n} hg_{\boldsymbol{\alpha}}.$$
\end{theorem}
\noindent\textbf{Remark 1.} By a theorem of Khintchine \cite{khin} and Theorem 2 in \cite{kem}, the set of such $\boldsymbol{\alpha}$ has full Lebesgue measure. 

\noindent\textbf{Remark 2.} It is possible to prove a similar conclusion
with a larger class of test functions h, using more complicated Fourier analysis techniques, along the lines of what was done in \cite{fz}, or by assuming the Riemann Hypothesis. 

\vspace{1em}
For the special case when $r=n$ and $Rank(M)=n$, \eqref{n-dim-alpha} will
always be satisfied as a consequence of Baker's Theorem \cite{Baker}, and we have
the following.
\begin{coro} \label{n-dim-as}
	If $h(\mathbf{x})$ satisfies condition (\ref{n-dim-condt}) and $\boldsymbol{\alpha}$ satisfies condition (\ref{n-dim-condition}) for $r=n$, we have
	$$\lim_{T\rightarrow\infty}\frac{1}{T}\left(\sum_{0<\gamma\leq T}h(\gamma\boldsymbol{\alpha})-N(T)\int_{\mathbb{T}^n} h\right)=\int_{\mathbb{T}^n} hg_{\boldsymbol{\alpha}}.$$
\end{coro}

In general, the conclusion of Theorem \ref{n-dim-aa} may fail if the 
Diophantine condition \eqref{n-dim-alpha} fails; for details see the Remarks at
the conclusion of Section \ref{section-2dim}.
For the case $n=2$, the behavior of $\sum h(\gamma \boldsymbol{\alpha})$ is determined
by the behavior of the convergents  $\frac{p_n}{q_n}$ of the continued fraction of $\xi:=\frac{\alpha_1}{\alpha_2}$.  Even if \eqref{n-dim-alpha} fails, which occurs when some $q_{n+1}$ are very large compared to $q_n$, we show that
 there are ``long'' intervals of $T$-values on which we may draw the conclusion of the 
previous theorem. 

\begin{theorem} \label{2-dim-range}
	Assume \eqref{n-dim-condt} and fix any $\epsilon>0$.  Let $U_{\boldsymbol{\alpha}}=\bigcup_{n=1}^{\infty}[q_n^{1+\epsilon}, e^{q_n^{B-\epsilon}}]$. Then,
	\begin{equation*}
	\lim_{\substack{T\in U_{\boldsymbol{\alpha}}\\ T\rightarrow\infty}}\frac{1}{T}\left(\sum_{0<\gamma\leq T}h(\alpha_1 \gamma, \alpha_2 \gamma)-N(T)\int_{\mathbb{T}^2} h\right)=\int_{\mathbb{T}^2}  hg_{\boldsymbol{\alpha}}.
	\end{equation*}
\end{theorem}

Finally, we show 
 some graphs comparing the limiting density function $g_{\boldsymbol{\alpha}}$ with numerical computations obtained from the first 100 million  zeros of the Riemann zeta-function (which were kindly provided by Tom\'as Oliveira e Silva, see \cite{Si}). 
Define
\begin{equation}\
	M(y_1, y_2; T)=M_{\boldsymbol{\alpha}}(y_1, y_2; T):=\frac{1}{T}\sum_{\substack{0<\gamma\leq T\\ \{\alpha_1\gamma\}<y_1, \{\alpha_2\gamma\}<y_2}} 1-y_1 y_2\frac{N(T)}{T}.\nonumber
\end{equation}

Take $T=42653549.761$, then $N(T)=10^8$. Denote $\Delta=\frac{1}{100}$. We partition $[0, 1)\times[0, 1)$ into $\Delta^{-2}$ small rectangles with side length $\Delta$. For some values of $\alpha_1$, $\alpha_2$, we plot the values of $DM:=\frac{1}{\Delta^2}(M(y_1+\Delta, y_2+\Delta)-M(y_1+\Delta, y_2)-M(y_1, y_2+\Delta)+M(y_1, y_2))$ for each small rectangle $[y_1, y_1+\Delta)\times[y_2, y_2+\Delta)$.

\begin{figure}[ht!]
  \centering
  \includegraphics[width=14cm]{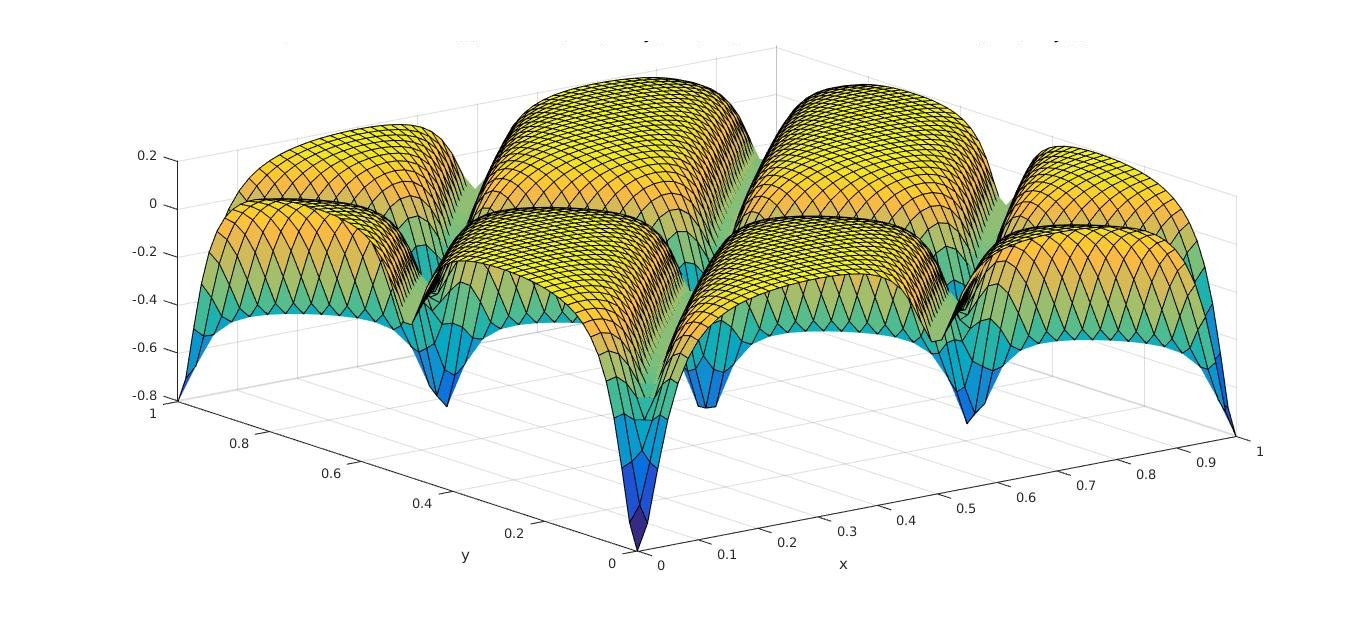}\\
  \caption{Graph of $g_{\boldsymbol{\alpha}}$}\label{g12}
\end{figure}

\begin{figure}[ht!]
  \centering
  \includegraphics[width=14cm]{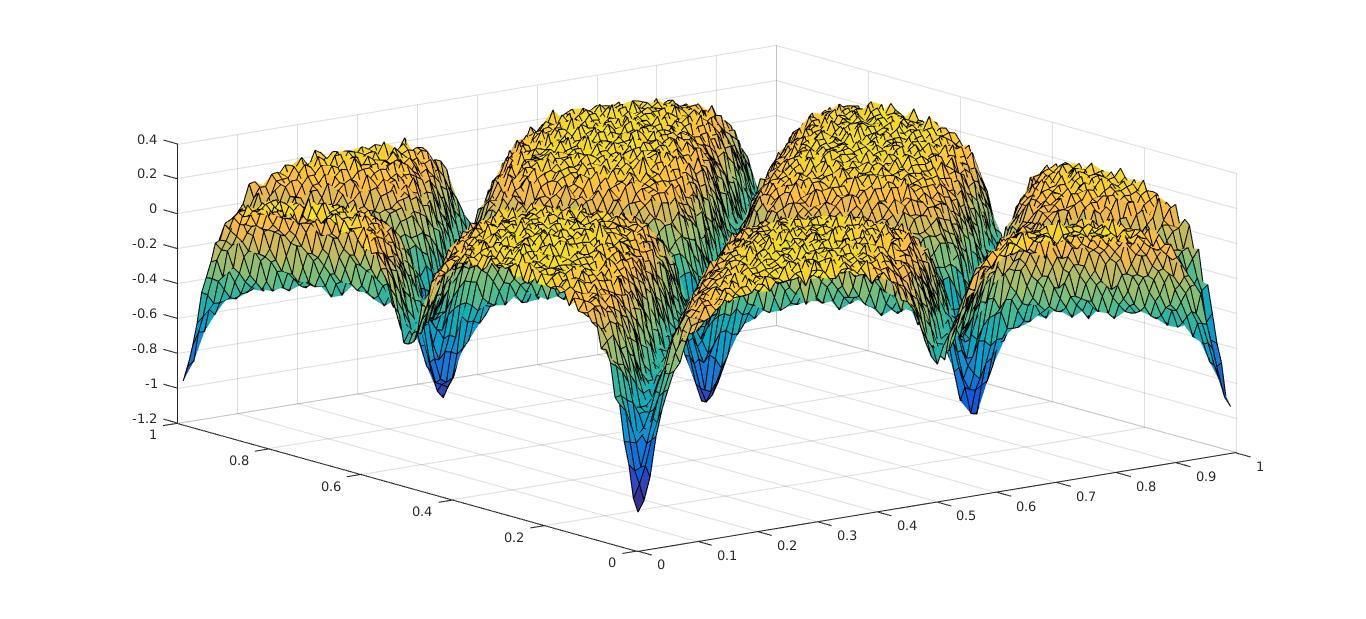}\\
  \caption{Values of $DM$}\label{22}
\end{figure}

\textbf{Example 1.} Let $\alpha_1$ and $\alpha_2$ satisfy 
\begin{align*}
\alpha_1+\alpha_2& =\dfrac{\log 2}{2\pi},\\
\alpha_1-\alpha_2& =\dfrac{\log 3}{4\pi}.
  \end{align*}
Figure $\ref{g12}$ shows the graph of $g_{\boldsymbol{\alpha}}$. Figure $\ref{22}$ shows the values of $DM$, the graph of which is a close approximation of  $g_{\boldsymbol{\alpha}}$.

\begin{figure}[ht!]
  \centering
  \includegraphics[width=14cm]{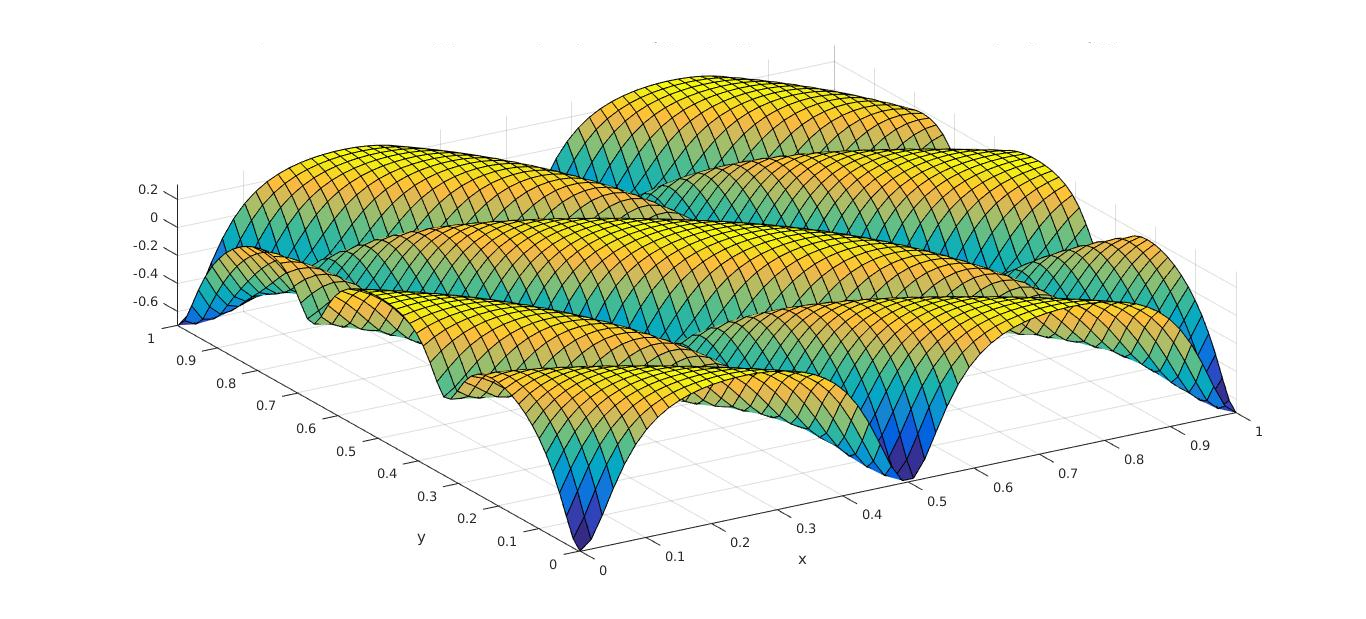}\\
  \caption{Graph of $g_{\boldsymbol{\alpha}}$}\label{23}
\end{figure}

\begin{figure}[ht!]
  \centering
  \includegraphics[width=14cm]{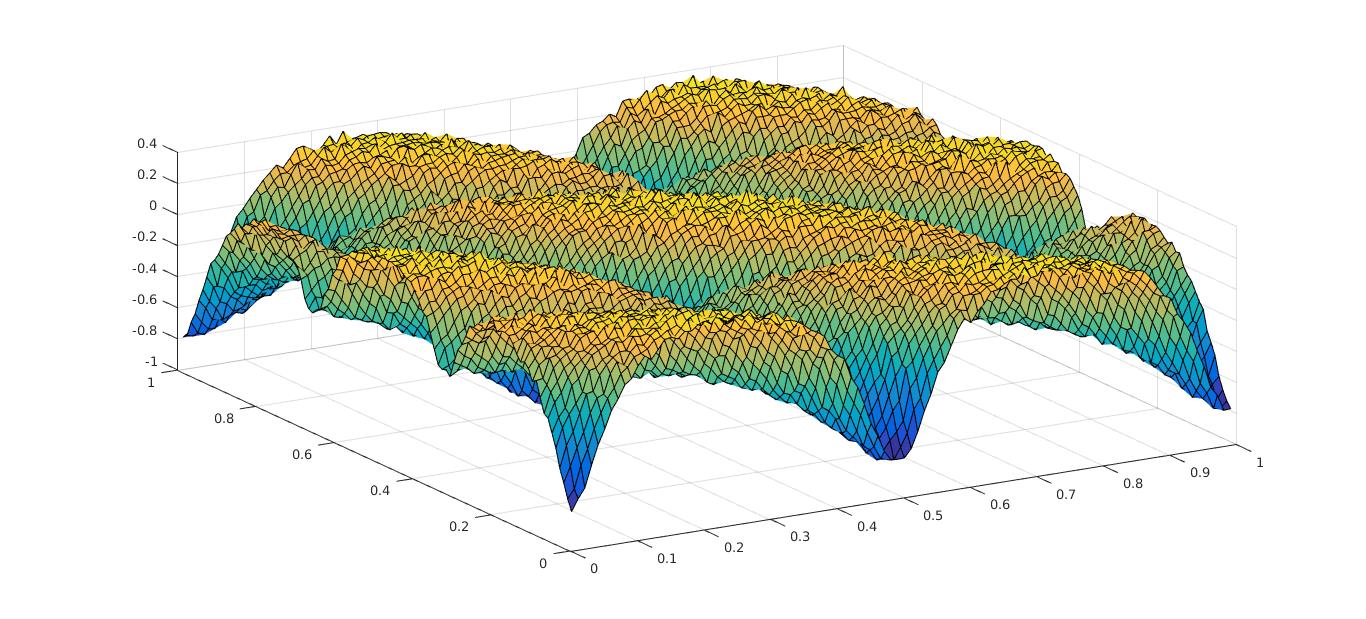}\\
  \caption{Values of $DM$}\label{24}
\end{figure}

\textbf{Example 2.} Figure $\ref{23}$ and Figure $\ref{24}$ show the analogous graphs for the case when $\alpha_1$ and $\alpha_2$ satisfy
\begin{align*}
  2\alpha_1+\alpha_2&=\dfrac{\log 5}{2\pi},\\
  2\alpha_1+3\alpha_2&=\dfrac{\log 7}{2\pi}.
\end{align*}

\FloatBarrier

\FloatBarrier

\FloatBarrier

\FloatBarrier

%
%

\section{Proof of Theorem \ref{n-dim-aa}}

%

We need Lemma 1 from \cite{fz}, which is an extension of a famous formula of Landau.  Let $\rho$ denote a generic nontrivial zero of the Riemann zeta-function, and denote $\gamma=\Im \rho$. 
\begin{lem}\label{lem-fz}
	Let $x, T>1$, and denote by $n_x$ the nearest prime power to $x$. Then
	$$\sum_{0<\gamma\leq T} x^{\rho}=-\frac{\Lambda(n_x)}{2\pi}\frac{e^{iT\log(x/n_x)}-1}{i\log(x/n_x)}+O\left(x\log^2(2xT)
	+\frac{\log 2T}{\log x}\right),$$
where if $x=n_x$ the first term is $-T\frac{\Lambda(n_x)}{2\pi}$.
\end{lem}

Next, we extract from the proof of (3.8) in \cite{fz} the following
estimate for difference between two sums over zeros, one with the 
real parts all normalized to be equal to 1/2.

\begin{lem}\label{lem-ford2}
	For $1<x\le \exp\{\frac{\log T}{50\log\log T} \}$, we have
	\begin{equation}\label{fz-relat}
	\sum_{0<\gamma\leq T}(x^{i\gamma}-x^{\rho-1/2})\ll \frac{T\log^2 x}{\log T}+\frac{T}{\log^{10} T}.
	\end{equation}
\end{lem}

\noindent \textbf{\textit{Proof of Theorem \ref{n-dim-aa}.}}

By Theorem 2 in \cite{kem}, we know that the set of $\boldsymbol{\alpha}$'s under the condition (\ref{n-dim-alpha}) has full Lebesgue measure. 

Write the Fourier expansion of $h(\mathbf{x})$ as
\begin{equation}\label{hidim-hfunc}
h(\mathbf{x})=\sum_{\mathbf{m}}c_{\mathbf{m}} e^{2\pi i(\mathbf{m}\cdot \boldsymbol{\alpha})}= \sum_{\|\mathbf{m}\|\leq J}c_{\mathbf{m}} e^{2\pi i (\mathbf{m}\cdot \boldsymbol{\alpha})}+O\left(\sum_{\|\mathbf{m}\|>J}|c_{\mathbf{m}}| \right),
\end{equation}
where $J\leq \log T$ is a parameter to be chosen later. 

By our assumption (\ref{n-dim-condt}), one can see that the error term above is $O(1/J^{B-n})$. Thus, by (\ref{hidim-hfunc}), we have
\begin{align*}
\sum_{0<\gamma\leq T} h(\gamma\boldsymbol{\alpha})&=\sum_{0<\gamma\leq T} \left( \sum_{\|\mathbf{m}\|\leq J}c_{\mathbf{m}} e^{2\pi i \gamma (\mathbf{m}\cdot \boldsymbol{\alpha})}\right)+O\left(\frac{N(T)}{J^{B-n}}\right)\\
&=N(T)\int h+2\Re \sum_{\substack{0<\|\mathbf{m}\|\leq J\\ \mathbf{m}\cdot \boldsymbol{\alpha}>0}}c_{\mathbf{m}}\sum_{0<\gamma\leq T}e^{2\pi i \gamma (\mathbf{m}\cdot \boldsymbol{\alpha})}+O\left(\frac{N(T)}{J^{B-n}}\right).
\end{align*}
The second equality is a consequence of the identity $c_{-\mathbf{m}}=\overline{c_{\mathbf{m}}}$. Let $x_{\mathbf{m}}=e^{2\pi (\mathbf{m}\cdot \boldsymbol{\alpha})}$. Hence, by \eqref{n-dim-condt} and Lemma \ref{lem-ford2}, we have for $T$ large and $J\ll \frac{\log T}{(\log\log T)^2}$ that
\begin{equation}\label{hidim-h}
\sum_{0<\gamma\leq T} h(\gamma\boldsymbol{\alpha})-N(T)\int h=2\Re \sum_{\substack{0<\|\mathbf{m}\|\leq J\\ \mathbf{m}\cdot \boldsymbol{\alpha}>0}}c_{\mathbf{m}}\sum_{0<\gamma\leq T}x_{\mathbf{m}}^{ \rho-1/2}+O\left(\frac{N(T)}{J^{B-n}}+\frac{T}{\log T}\right).
\end{equation}
Write
$$I=\sum_{\substack{0<\|\mathbf{m}\|\leq J\\ \mathbf{m}\cdot \boldsymbol{\alpha}>0}}c_{\mathbf{m}}\sum_{0<\gamma\leq T}x_{\mathbf{m}}^{ \rho-1/2}.$$
Applying Lemma \ref{lem-fz}, we get 
\begin{eqnarray}\label{n-dim-lem}
I&=&-\frac{1}{2\pi}\sum_{\substack{\mathbf{m}\cdot\boldsymbol{\alpha}>0\\ 0<\|\mathbf{m}\|\leq J}}\frac{c_{\mathbf{m}}\Lambda(n_{x_{\mathbf{m}}})}{\sqrt{x_{\mathbf{m}}}}   \frac{\sin(T\log{\frac{x_{\mathbf{m}}}{n_{x_{\mathbf{m}}}}})}{\log{\frac{x_{\mathbf{m}}}{n_{x_{\mathbf{m}}}}}}      \nonumber\\
&& + O\left(\sum_{\substack{\mathbf{m}\cdot\boldsymbol{\alpha}>0\\ 0<\|\mathbf{m}\|\leq J}} |c_{\mathbf{m}}| \left(\sqrt{x_{\mathbf{m}}}\log^2(2x_{\mathbf{m}}T)+\frac{\log 2T}{\sqrt{x_{\mathbf{m}}}\log x_{\mathbf{m}}}  \right) \right)\nonumber\\
&=&-\frac{1}{2\pi}\sum_{\substack{\mathbf{m}\cdot\boldsymbol{\alpha}>0\\ 0<\|\mathbf{m}\|\leq J}}\frac{c_{\mathbf{m}}\Lambda(n_{x_{\mathbf{m}}})}{\sqrt{x_{\mathbf{m}}}}   \frac{\sin(T\log{\frac{x_{\mathbf{m}}}{n_{x_{\mathbf{m}}}}})}{\log{\frac{x_{\mathbf{m}}}{n_{x_{\mathbf{m}}}}}} +O\left(e^{O(J)}\log^2 T\right).
\end{eqnarray}
Here we used our assumption (\ref{n-dim-alpha}) that
$$\log x_{\mathbf{m}}=2\pi (\mathbf{m}\cdot\boldsymbol{\alpha}) \geq \frac{2\pi C}{e^{J}}.$$

Now let $J=\sqrt{\log T}$. By (\ref{hidim-h}), (\ref{n-dim-lem}) and $B>n+2$, we get
\begin{equation}\label{hidim-t}
\sum_{0<\gamma\leq T} h(\gamma\boldsymbol{\alpha})-N(T)\int h =-\frac{T}{\pi}\Re\sum_{\substack{\mathbf{m}\cdot\boldsymbol{\alpha}>0\\ 0<\|\mathbf{m}\|\leq J}}\frac{c_{\mathbf{m}}\Lambda(n_{x_{\mathbf{m}}})}{\sqrt{x_{\mathbf{m}}}}   \frac{\sin(T\log{\frac{x_{\mathbf{m}}}{n_{x_{\mathbf{m}}}}})}{T\log{\frac{x_{\mathbf{m}}}{n_{x_{\mathbf{m}}}}}}+o(T).
\end{equation}
Since $\frac{\Lambda(n_{x_{\mathbf{m}}})}{\sqrt{x_{\mathbf{m}}}}$ is bounded, the first sum is absolutely and uniformly convergent in $T$ under the assumption $(\ref{n-dim-condt})$. And the terms with $n_{x_{\mathbf{m}}}\neq x_{\mathbf{m}}$ tend to zero as $T\rightarrow\infty$. Thus, by (\ref{hidim-t}), as $T\rightarrow\infty$,
\begin{equation}\label{n-dim-conclu}
\frac{1}{T}\left( \sum_{0<\gamma\leq T} h(\gamma\boldsymbol{\alpha})-N(T)\int h\right)
=-\frac{1}{\pi}\Re \left\{\sum_{\substack{\mathbf{m}\cdot\boldsymbol{\alpha}>0\\0<\|\mathbf{m}\|\leq J}}\frac{c_{\mathbf{m}}\Lambda({x_{\mathbf{m}}})}{\sqrt{x_{\mathbf{m}}}}         \right\}+o(1).
\end{equation}
Since $\Lambda({x_{\mathbf{m}}})=0$ unless $x_{\mathbf{m}}$ is prime power, by (\ref{n-dim-g2}), when $g_{\boldsymbol{\alpha}}=0$, the sum in the right hand side of (\ref{n-dim-conclu}) is also $0$. If $\boldsymbol{\alpha}$ satisfies (\ref{n-dim-condition}) with $M$ being maximal (maximal $r$), then the only terms with 
$n_{x_{\mathbf{m}}} =  x_{\mathbf{m}}$ (that is, with $x_{\mathbf{m}}$ being a prime power)
are those with $\mathbf{m}$ being a multiple of some row of $M$.
Hence, recalling the definition \eqref{n-dim-g1} of $g_{\boldsymbol{\alpha}}$, as
$T\rightarrow \infty$, the right hand side of (\ref{n-dim-conclu}) becomes 
$$-\frac{1}{\pi}\Re\left\{\sum_{j=1}^r \sum_{k\geq 1}(\log p_j) p_j^{-\frac{a_j k}{2}} c_{k{q}_j \boldsymbol{m}}\right\},$$
which equals $\int_{\mathbb{T}^n} hg_{\boldsymbol{\alpha}}$ by (\ref{n-dim-g1}). 
So in all cases, we have
$$\lim_{T\rightarrow\infty}\frac{1}{T}\left(\sum_{0<\gamma\leq T}h(\gamma\boldsymbol{\alpha})-N(T)\int_{\mathbb{T}^n} h\right)=\int_{\mathbb{T}^n} hg_{\boldsymbol{\alpha}}.$$

Therefore,  we get the conclusion of this theorem. 

\vspace{1em}

\noindent\textbf{\textit{Proof of Corollary \ref{n-dim-as}.}}

For the special form of $\boldsymbol{\alpha}$ when $r=n$,  we can solve the above linear equations for $\boldsymbol{\alpha}$, and find that each $2\pi\alpha_i$ is a linear combination of these $\log p_i$'s with rational coefficients. Hence, by Baker's theorem \cite{Baker}, for such $\boldsymbol{\alpha}$, there exist constants  $D=D(\boldsymbol{\alpha})$  and $\mu$, such that, for all $\mathbf{m}\in \mathbb{Z}^n$, 
\begin{equation}\label{n-dim-Baker}
|\mathbf{m}\cdot \boldsymbol{\alpha}|\geq \frac{D}{(\|\mathbf{m}\|+1)^\mu}.
\end{equation}
Thus, such $\boldsymbol{\alpha}$ satisfies condition (\ref{n-dim-alpha}) for some constant $C$. 
The conclusion follows from Theorem \ref{n-dim-aa}.

\section{Proof of Theorem \ref{2-dim-range}}\label{section-2dim}
We may assume that $\alpha_1/\alpha_2$ is positive and irrational, else the problem reduces to the $n=1$ case.

Assume a function $h(x,y)$ is defined on the two dimensional torus. By (\ref{hidim-h}), we have
\begin{equation}\label{star-h2}
 \sum_{0<\gamma\leq T}h(\alpha_1 \gamma, \alpha_2 \gamma)-N(T)\int h = 2\Re \sum_{\substack{|m|,|l|\le J \\ m\alpha_1+l\alpha_2>0}} c_{m,l}
 \sum_{0<\gamma\leq T} x_{m,l}^{\rho-1/2} + O\left(\frac{N(T)}{J^{B-2}}+\frac{T}{\log T}\right),
\end{equation}
where $J\leq \log T$ is a parameter to be chosen later. 

We now break the double sum into two pieces, over those pairs $(m,l)$ corresponding to ``large'' $x_{m,l}$ and those
pairs $(m,l)$ corresponding to ``small'' $x_{m,l}$.
Fix an arbitrary positive constant $C$ and denote
\begin{align*}
  E_J&:=\left\{(m, l)\in \mathbb{Z}^2:0<\max\{|m|, |l|\}\leq J, m\alpha_1+l\alpha_2>\min\(\frac{C}{e^{\max\{|m|, |l|\}}}, \frac{|\alpha_2|}{2m}, \frac1{4\pi}\)\right\}, \\
  F_J&:=\left\{(m, l)\in \mathbb{Z}^2:0<\max\{|m|, |l|\}\leq J, 0 < m\alpha_1+l\alpha_2 \le \min\(\frac{C}{e^{\max\{|m|, |l|\}}}, \frac{|\alpha_2|}{2m}, \frac1{4\pi}\)\right\}.
\end{align*}

Write
\[
I_E=\sum_{(m,l)\in E_J}c_{m,l}\sum_{0<\gamma\leq T} x_{m,l}^{\rho-1/2}, \qquad 
I_F=\sum_{(m,l)\in F_J}c_{m,l}\sum_{0<\gamma\leq T} x_{m,l}^{\rho-1/2}.
\]
Applying Lemma \ref{lem-fz} and using $B>4$,  we get
\begin{eqnarray}\label{est-I1}
  I_E&=&-\frac{1}{2\pi} \sum_{(m,l)\in E_J}\frac{c_{m,l}\Lambda(n_{x_{m,l}})}{\sqrt{x_{m,l}}}\frac{\sin(T\log\frac{x_{m,l}}{n_{x_{m,l}}})}{\log\frac{x_{m,l}}{n_{x_{m,l}}}}\nonumber\\
  &&+O\left(\sum_{(m, l)\in E_J}|c_{m,l}|\left(\sqrt{x_{m,l}}\log^2(2x_{m,l}T)+\frac{\log 2T}{\sqrt{x_{m,l}}\log x_{m,l}}\right)\right)\nonumber\\
  &=&-\frac{1}{2\pi} \sum_{(m,l)\in E_J}\frac{c_{m,l}\Lambda(n_{x_{m,l}})}{\sqrt{x_{m,l}}}\frac{\sin(T\log\frac{x_{m,l}}{n_{x_{m,l}}})}{\log\frac{x_{m,l}}{n_{x_{m,l}}}}
  +O\left(e^{O(J)}\log^2 T\right).
\end{eqnarray}
Here we have used the fact that for $(m,\ell)\in E_J$,
\[
 \log x_{m,\ell} = 2\pi (m\alpha_1+l\alpha_2) \gg e^{-J}.
\]

By standard facts on Diophantine approximation, we know that if $(m,l)\in F_J$ then $$\left| \frac{\alpha_1}{\alpha_2}+\frac{l}{m}\right| < \frac{1}{2m^2},$$
and hence all elements of $F_J$ are of the form $(q_n, -p_n)$, where $p_n/q_n$ is a convergent 
of the continued fraction of $\frac{\alpha_1}{\alpha_2}$.
 As $q_n$ grows at least exponentially and $q_n\le J$ for $(q_n,-p_n)\in F_J$,
there are $O(\log J)=O(\log\log T)$ elements in $F_J$.  Let $n^*$ denote the largest $n$
such that $(q_n,-p_n)\in F_J$.

By elementary properties of  continued fractions (Chapter X, \cite{hw}, or Corollary 1.4, \cite{bu}), we have
\begin{equation}\label{cont-fracp}
  \frac{\alpha_2}{q_j+q_{j+1}}<|q_j \alpha_1-p_j \alpha_2|<\frac{\alpha_2}{q_{j+1}}.
\end{equation}
Thus, for all $j$ except possibly $j=n^*$, we have $q_{j+1} \le J$ and
from Lemma \ref{lem-fz} we get (writing $x_j=x_{q_j,-p_j}$) 
\begin{align}\label{est-I2}
\sum_{j\ne n^*}c_{q_j,-p_j}\sum_{0<\gamma\leq T} x_j^{\rho-1/2}&=-\frac{1}{2\pi} 
\sum_{j\ne n^*}\frac{c_{q_j,-p_j}\Lambda(n_{x_j})}{\sqrt{x_j}}\frac{\sin(T\log\frac{x_j}{n_{x_j}})}{\log\frac{x_j}{n_{x_j}}} 
+O(\log^3 T) \nonumber\\
&=O(\log^3 T),
\end{align}
where we have used that $n_{x_j}=2$. We adopt the convention that the term $j=n^*$ is included in the sum above 
if $q_{n^*+1} \le J$. 
In fact, the error term from Lemma \ref{lem-fz} is acceptable, namely $o(T)$, if $q_{n^*+1}=o(T/\log T)$ as well.

Now let $J=\sqrt{\log T}$ and assume that $(\log T)^{1/B}\log\log T< q_{n^*} \le J$. 
We revert to the sum over $x_{q_{n^*},-p_{n^*}}^{i\gamma}$.  By \eqref{fz-relat} and the the trivial
bound for $\sum_{0 <\gamma\le T} x_{q_{n^*},-p_{n^*}}^{i\gamma}$, we have
$$c_{q_{n^*}, -p_{n^*}}\sum_{0<\gamma\leq T} x_{q_{n^*},-p_{n^*}}^{\rho-1/2} \ll 
\frac{T\log T}{q_{n^*}^B} =o(T).$$
Combined with \eqref{est-I1} and \eqref{est-I2}, we see that if (a) $F_J$ is empty or if (b)
there is some convergent of $\alpha_1/\alpha_2$ with
$q_n \in ((\log T)^{1/B}\log\log T, o(T/\log T)]$, then
\begin{eqnarray}\label{star-I}
  \sum_{0<\gamma\leq T}h(\alpha_1 \gamma, \alpha_2 \gamma)&=&-\frac{T}{\pi}\Re \sum_{\substack{m\alpha_1+l\alpha_2>0\\0<\max\{|m|, |l|\}\leq J}}\frac{c_{m,l}\Lambda(n_{x_{m,l}})}{\sqrt{x_{m,l}}}\frac{\sin(T\log\frac{x_{m,l}}{n_{x_{m,l}}})}{T\log\frac{x_{m,l}}{n_{x_{m,l}}}}\nonumber\\
  &&+N(T) \int h + o(T).
\end{eqnarray}
Note that (b) is satisfied  if $T$ is large enough and
$T\in \bigcup_{n=1}^{\infty}[q_n^{1+\epsilon}, e^{q_n^{B-\epsilon}}]$.

Then, similar to the proof of (\ref{n-dim-conclu}) in the proof of Theorem \ref{n-dim-aa}, we deduce that
\begin{equation*}
  \lim_{\substack{T\in U_{\boldsymbol{\alpha}}\\ T\rightarrow\infty}}\frac{1}{T}\left(\sum_{0<\gamma\leq T}h(\alpha_1 \gamma, \alpha_2 \gamma)-N(T)\int h\right)=\int hg_{\boldsymbol{\alpha}}.
\end{equation*}
Thus, Theorem \ref{2-dim-range} follows.

\medskip

\textbf{Remarks.}  Suppose that there is a very long gap between convergents,
say $q_{n^*}$ is small and $q_{n^*+1}$ is very large.  Put $x=x_{q_{n^*},-p_{n^*}}$ and assume that
$\log x \ll \frac{\log T}{T}$.  Using the known estimate $N(T)=\frac{T}{2\pi}\log(\frac{T}{2\pi e})+O(\log T)$,
we get
\begin{align*}
 \sum_{0<\gamma\le T} x^{i\gamma} &=\int_0^T x^{it} \log \left(\frac{t}{2\pi}\right)\, dt + O(\log^2 T)\\
 &=T\log(T/2\pi) \frac{e^{iT\log x}-1}{iT\log x} + O(T).
\end{align*}
In the range $T\ll 1/\log x$, there is thus a term in the sum \eqref{star-h2} of order $N(T)|c_{q_{n^*},-p_{n^*}}|$,
which may be larger than $T$.

\vspace{1em}
\textbf{Acknowledgement.} Research of the first author and second author is partially
supported by NSF grants DMS-1201442 and DMS-1501982. The authors would like to thank the helpful comments of the referee.

{\footnotesize
Department of Mathematics, University of Illinois at Urbana-Champaign, 1409 West Green Street, Urbana, IL 61801, USA

\textit{E-mail}: KF: ford@math.uiuc.edu, XM: xmeng13@illinois.edu, 
\newline 

Department of Mathematics, University of Illinois at Urbana-Champaign, 1409 West Green Street,
Urbana, IL 61801, USA and Institute of Mathematics of the Romanian
Academy, P.O. Box 1-764, Bucharest RO-70700, Romania.

\textit{E-mail}: AZ: zaharesc@math.uiuc.edu}

\end{document}